\documentclass[11pt,twoside]{article}

\usepackage{a4wide}
\usepackage{amsfonts}
\usepackage{amssymb}
\usepackage{amsmath}
\usepackage{graphicx}

\newcommand{\ignore}[1]{}

\def\@begintheorem#1#2{\par\bgroup{\sc #1\ #2. }\it\ignorespaces}
\def\@opargbegintheorem#1#2#3{\par\bgroup{\sc #1\ #2\ (#3). } \it\ignorespaces}
\def\@endtheorem{\egroup}
\newtheorem{theorem}{Theorem}[section]
\newtheorem{corollary}[theorem]{Corollary}
\newtheorem{lemma}[theorem]{Lemma}

\newtheorem{example}[theorem]{Example}
\newtheorem{proposition}[theorem]{Proposition}
\newtheorem{definition}[theorem]{Definition}
\newcommand{\bt}[1]{\begin{theorem}\label{#1}}
\newcommand{\bc}[1]{\begin{corollary}\label{#1}}
\newcommand{\bl}[1]{\begin{lemma}\label{#1}}
\newcommand{\be}[1]{\begin{example}\label{#1}}
\newcommand{\bp}[1]{\begin{proposition}\label{#1}}
\newcommand{\ba}[1]{\begin{algorithm}\rm\label{#1}}
\newcommand{\bd}[1]{\begin{definition}\rm\label{#1}}{\normalsize }
\newcommand{\bpr}{\noindent {\em Proof. }}
\newcommand{\et}{\end{theorem}}
\newcommand{\ec}{\end{corollary}}
\newcommand{\el}{\end{lemma}}
\newcommand{\ee}{\end{example}}
\newcommand{\ep}{\end{proposition}}
\newcommand{\ed}{\end{definition}}
\newcommand{\epr}{{\ \vbox{\hrule\hbox{%
\vrule height1.3ex\hskip0.8ex\vrule}\hrule}}\\\par}

\def\Z{\mathbb{Z}}

\begin{document}

\title{\bf On Line Sum Optimization}

\author{
Shmuel Onn
\thanks{\small Technion - Israel Institute of Technology.
Email: onn@ie.technion.ac.il}
}
\date{}

\maketitle

\begin{abstract}
We show that the {\em column sum optimization problem}, of finding a $(0,1)$-matrix
with prescribed row sums which minimizes the sum of evaluations of given functions
at its column sums, can be solved in polynomial time, either when all functions
are the same or when all row sums are bounded by any constant.
We conjecture that the more general {\em line sum optimization problem},
of finding a matrix minimizing the sum of given functions evaluated at
its row sums and column sums, can also be solved in polynomial time.


\end{abstract}

\section{Introduction}

For a positive integer $m$ let $[m]:=\{1,2,\dots,m\}$.
For an $m\times n$ matrix $A$ let $r_i(A):=\sum_{j=1}^n A_{i,j}$ for $i\in[m]$
and $c_j(A):=\sum_{i=1}^m A_{i,j}$ for $j\in[n]$ be its row sums and column sums.

We consider here the following algorithmic problem.

\vskip.2cm\noindent{\bf Column Sum Optimization.}
Given positive integers $m,n$ and $r_1,\dots,r_m\leq n$, and functions
$f_j:\{0,1,\dots,m\}\rightarrow\Z$ for $j=1,\dots,n$, find an $m\times n$ matrix $A$,
where each entry is equal to $0$ or $1$, and with row sums $r_i(A)=r_i$ for $i=1,\dots,m$,
which minimizes $\sum_{j=1}^n f_j(c_j(A))$.

\vskip.2cm
For instance, for $m=n=4$, row sum tuple $(r_1,r_2,r_3,r_4)=(3,3,2,2)$, and identical functions
$f_1=f_2=f_3=f_4=f$ with $f(x)=(x-1)^2(x-3)^2$, an optimal solution is the
following matrix, with column sum tuple $(c_1,c_2,c_3,c_4)=(3,3,3,1)$ and objective value $0$,
$$A\ = \left(\begin{array}{cccc}
1&1&1&0\\
1&1&1&0\\
1&1&0&0\\
0&0&1&1\\
\end{array}
\right)\ .$$

\vskip.2cm
We use a characterization of $(0,1)$-matrices by Ryser \cite{Rys} and prove the following theorem.
\bt{main}
The column sum optimization problem can be solved in polynomial time when:
\begin{enumerate}
\item
either all functions are the same, that is, $f_1=\cdots=f_n=f$ for some $f$;
\item
or all row sums are bounded by some fixed integer $b$, that is, $r_1,\dots,r_m\leq b$.
\end{enumerate}
\et

Note that even when the row sums are bounded, the number of possible $(0,1)$-matrices can be
${n\choose b}^m$ which is exponential in the data, hence exhaustive search is not polynomial.

\vskip.2cm
The column sum optimization problem is a special case of the following more general
{\em line sum optimization problem}: given $m,n$ and $e\leq mn$, and functions
$g_i:\{0,1,\dots,n\}\rightarrow\Z$ for $i\in[m]$ and $f_j:\{0,1,\dots,m\}\rightarrow\Z$
for $j\in[n]$, find $A\in\{0,1\}^{m\times n}$ with $e$ nonzero entries minimizing
$\sum_{i=1}^m g_i(r_i(A))+\sum_{j=1}^n f_j(c_j(A))$. Indeed, the column sum problem
reduces to the line sum problem by setting $e:=\sum_{i=1}^mr_i$ and $g_i(x):=a(x-r_i)^2$
for $i\in[m]$ and sufficiently large positive integer $a$.
However, we do not know what is the complexity of this more general problem,
and not even that of the column sum problem for arbitrary row sums and functions.

\vskip.2cm
The line sum problem is in turn a special case of the following
{\em degree sequence optimization problem}: given a graph $H=(V,E)$, integer $e\leq|E|$,
and functions $f_v:\{0,1,\dots,d_v(H)\}\rightarrow\Z$ for $v\in V$, where $d_v(H)$
is the degree of $v$ in $H$, find a subgraph $G=(V,F)\subseteq H$
with $e$ edges minimizing $\sum_{v\in V}f_v(d_v(G))$. Indeed, identifying matrices
$A\in\{0,1\}^{m\times n}$ with graphs $G=(V,F)$ where
$V=\{u_1,\dots,u_m\}\uplus\{w_1,\dots,w_n\}$ and $F=\{\{u_i,w_j\}\,:\,A_{i,j}=1\}$,
the line sum problem reduces to the degree sequence problem with $H=K_{m,n}$ the
complete bipartite graph. In the case of $H=K_n$ the complete graph and all functions the same,
$f_v=f$ for all $v\in V$, the problem was recently shown in \cite{DLMO} to be polynomial time
solvable, using the characterization of degree sequences by Erd\H{o}s and Gallai \cite{EG}.
For general graphs $H$, the problem was shown in \cite{AS} to be NP-hard already when
$f_v(x)=-x^2$ for all $v\in V$, but is polynomial time solvable if the functions are
convex \cite{AS,DO}. If $m$ or $n$ are fixed, the problem for any $H\subseteq K_{m,n}$, and hence in particular the line sum problem, for any functions, is
polynomial time solvable \cite{DO}. We conjecture that the line sum problem, which is the degree sequence
problem over $H=K_{m,n}$, as well as the degree sequence problem over $H=K_n$, are
polynomial time solvable for arbitrary functions at the vertices, not necessarily identical.

\section{Proof}

Let the positive integers $m,n$ and $r_1,\dots,r_m\leq n$ be given.
Define the tuple $r:=(r_1,\dots,r_m)$ and its {\em conjugate} tuple
$d:=(d_1,\dots,d_n)$ where $d_j:=|\{i\,:\,r_i\geq j\}|$ for $j\in[n]$. Note that
$\sum_{j=1}^n d_j=\sum_{i=1}^m r_i\leq mn$ and $d$ is nonincreasing, that is,
$d_1\geq\cdots\geq d_n$. A nonincreasing tuple $c=(c_1,\dots,c_n)$ is {\em majorized} by $d$ if
$\sum_{j=1}^k c_j\leq\sum_{j=1}^k d_j$ for $k\in[n]$ and $\sum_{j=1}^n c_j=\sum_{j=1}^n d_j$.
(See \cite{MOA} for more details on the theory of majorization and its applications.)

\vskip.2cm
We use a characterization by Ryser \cite{Rys} of $(0,1)$-matrices with prescribed line sums.
\bp{Ryser}{\bf \cite{Rys}}
There is a matrix $A\in\{0,1\}^{m\times n}$ with row sum and column sum tuples $r,c$
if and only if, permuting $c$ to be nonincreasing, it is majorized by the conjugate tuple $d$.
\ep

For instance, if $m=n=4$ and $r=(3,3,2,2)$, then $d=(4,4,2,0)$, and $c=(3,3,3,1)$ is majorized by
$d$ so there exists $A\in\{0,1\}^{4\times 4}$ with row sum tuple $r$ and column sum tuple $c$.

\vskip.2cm
Next we note that a matrix with given line sums $r$ and $c$ can be efficiently obtained.
\bl{flow}
Let $r=(r_1,\dots,r_m)$ and let $c=(c_1,\dots,c_n)$ be such that,
permuting $c$ to be nonincreasing, it is majorized by the tuple
$d=(d_1,\dots,d_n)$ conjugate to $r$. Then a matrix $A$ in $\{0,1\}^{m\times n}$ with
row sum tuple $r$ and column sum tuple $c$ is computable in polynomial time.
\el

\newpage

\bpr
The problem is solvable either by the efficient simple Gale-Ryser algorithm,
see \cite[Chapter 3]{Bru}, or by network flows as follows.
Define a directed graph with capacities on the edges as follows.
There are vertices $s,t$, $u_1,\dots,u_m$, and $w_1,\dots,w_n$. There are edges $[s,u_i]$ for
$i\in[m]$ with capacity $r_i$, edges $[u_i,w_j]$ for $i\in[m]$ and $j\in[n]$ with capacity $1$,
and edges $[w_j,t]$ for $j\in[n]$ with capacity $c_j$. Then, as is well known, a maximum
nonnegative integer flow from $s$ to $t$ can be computed in polynomial time, see e.g. \cite{Sch}.
Then $A$ is read off from the maximum flow by taking $A_{i,j}$
to be the flow on edge $[u_i,w_j]$ for all $i\in[m]$ and $j\in[n]$.
\epr
We are now in position to prove our theorem.

\vskip.2cm
{\em Part 1.} Since all functions at the column sums are now the same, that is, $f_1=\cdots=f_n=f$
for some $f$, the objective value $\sum_{j=1}^n f(c_j(A))$ is invariant under column
permutations, and so without loss of generality we may and do search for a matrix $A$
whose column sum tuple $c$ is nonincreasing. (This is not true, however, for general
objective $\sum_{j=1}^n f_j(c_j(A))$ with distinct $f_j$, and this is precisely
the reason we do not know how to solve the general problem.)

\vskip.2cm
We now reduce the column sum problem to the problem of finding a shortest directed path
in a directed graph $D$ where each edge has a length. We construct $D$ as follows.

There are two special vertices $u$ and $v$. In addition, there are vertices labeled
by triplets of integers $(k,c_k,s_k)$ where $1\leq k\leq n$, $0\leq c_k\leq d_1\leq m$,
and $0\leq s_k\leq\sum_{j=1}^nd_j\leq mn$.

There are edges $[u,(1,c_1,s_1)]$ for $0\leq c_1\leq d_1$ and $s_1=c_1$ of length $f(c_1)$.
For $2\leq k\leq n$ there are edges $[(k-1,c_{k-1},s_{k-1}),(k,c_k,s_k)]$
for $c_k\leq c_{k-1}$, $s_k=s_{k-1}+c_k$, and $s_k\leq \sum_{j=1}^k d_j$, of length $f(c_k)$.
Finally, there are edges $[(n,c_n,s_n),v]$ for $s_n=\sum_{j=1}^n d_j$ of length $0$.

\vskip.2cm
Consider any directed path from $u$ to $v$ in $D$, which by the definition of $D$ looks like
$$u\longrightarrow(1,c_1,s_1)\longrightarrow\cdots\cdots
\longrightarrow(n,c_n,s_n)\longrightarrow v\ .$$

Intuitively, traveling on the path, being at vertex $(k,c_k,s_k)$ indicates that the prefix
$(c_1,\dots,c_k)$ of the sequence had been determined, and it is nonincreasing and sums to $s_k$.
Formally, by definition of the edges in $D$, we have $c_1\geq\cdots\geq c_n$ hence the
corresponding tuple $c=(c_1,\dots,c_n)$ is nonincreasing.
We also have that $\sum_{j=1}^k c_j=s_k\leq \sum_{j=1}^k d_j$ for $k\in[n]$ and
$\sum_{j=1}^n c_j=s_n=\sum_{j=1}^n d_j$, and therefore $c$ is majorized by the conjugate
tuple $d$. So by Proposition \ref{Ryser}, there is a matrix $A\in\{0,1\}^{m\times n}$
with row sum tuple $r$ and column sum tuple $c$. By the definition of the lengths,
the length of this directed path is $\sum_{j=1}^n f(c_j)=\sum_{j=1}^n f(c_j(A))$.

Conversely, for any matrix $A\in\{0,1\}^{m\times n}$ with row sum tuple $r$ and nonincreasing
column sum tuple $c=(c_1,\dots,c_n)$, the tuple $c$ is majorized by $d=(d_1,\dots,d_n)$
by Proposition \ref{Ryser}, hence the following is a directed path from
$u$ to $v$ in $D$ of length $\sum_{j=1}^n f(c_j)=\sum_{j=1}^n f(c_j(A))$,
$$u\longrightarrow(1,c_1,c_1)\longrightarrow(1,c_2,c_1+c_2)\longrightarrow\cdots\cdots
\longrightarrow(n,c_n,\sum_{j=1}^n c_j)\longrightarrow v\ .$$

So the column sum problem indeed reduces to the shortest directed path problem.

\vskip.2cm
Now, the number of vertices of $D$ is bounded by $n(m+1)(mn+1)+2=O(m^2n^2)$ and hence is
polynomial in the input. Therefore, as is well known, a shortest directed path from $u$ to $v$
in $D$ can then be obtained in polynomial time, see e.g. \cite{Sch}. Then an optimal
nonincreasing tuple $c=(c_1,\dots,c_n)$ of column sums can be read off from this optimal
directed path as explained above. Finally, an optimal matrix $A\in\{0,1\}^{m\times n}$,
with row sum tuple $r=(r_1,\dots,r_m)$
and column sum tuple $c=(c_1,\dots,c_n)$,
can be computed in polynomial time by Lemma \ref{flow}.
\epr

\vskip.2cm
{\em Part 2.} Since all row sums are now bounded above by the constant $b$, by definition
of the tuple $d=(d_1,\dots,d_n)$ conjugate to $r$ we now have that $d_j=|\{i\,:\,r_i\geq j\}|=0$
for all $j>b$. This implies that $\sum_{j=1}^b d_j=\sum_{j=1}^n d_j$ and hence a
nonincreasing tuple $c=(c_1,\dots,c_n)$ is majorized by $d=(d_1,\dots,d_n)$
already if $\sum_{j=1}^i c_j\leq\sum_{j=1}^i d_j$ for $i=1,\dots,b$ and
$\sum_{j=1}^n c_j=\sum_{j=1}^n d_j$, since for $i>b$ we then automatically
have $\sum_{j=1}^i c_j\leq\sum_{j=1}^n c_j=\sum_{j=1}^n d_j=\sum_{j=1}^i d_j$.

Note that, unlike the assumption in Part 1, the functions $f_j$ can now be distinct, so we cannot
assume that the matrix we are looking for has nonincreasing column sums.

Again we reduce the problem to the problem of finding a shortest directed path
in a directed graph $D$ where each edge has a length. We now construct $D$ as follows.

There are two vertices $u$ and $v$. There are also vertices labeled by $(b+3)$-tuples
of integers $(k,c_k,c^k_1,\dots,c^k_b,s_k)$ where $1\leq k\leq n$, $0\leq c_k\leq d_1\leq m$,
$m\geq d_1\geq c^k_1\geq\cdots\geq c^k_b\geq 0$, and $0\leq s_k\leq\sum_{j=1}^n d_j\leq mn$.
There are edges $[u,(1,c_1,c^1_1,\dots,c^1_b,s_1)]$ for $0\leq c_1\leq d_1$,
$c^1_1=c_1,c^1_2=\cdots=c^1_b=0$, and $s_1=c_1$, of length $f_1(c_1)$.
For $2\leq k\leq n$ there are edges
$$[(k-1,c_{k-1},c^{k-1}_1,\dots,c^{k-1}_b,s_{k-1}),(k,c_k,c^k_1,\dots,c^k_b,s_k)]\ ,$$
where $c^k_1\geq\cdots\geq c^k_b$ are the $b$ largest values among the $b+1$ values
$c^{k-1}_1,\dots,c^{k-1}_b,c_k$, where $s_k=s_{k-1}+c_k$,
and where $\sum_{j=1}^i c^k_j\leq\sum_{j=1}^i d_j$ for $i=1,\dots,b$,
of length $f_k(c_k)$. Finally, there are edges $[(n,c_n,c^n_1,\dots,c^n_b,s_n),v]$
for $s_n=\sum_{j=1}^n d_j$ of length $0$.

\vskip.2cm
Consider any directed path from $u$ to $v$ in $D$, which by the definition of $D$ looks like
$$u\longrightarrow(1,c_1,c^1_1,\dots,c^1_b,s_1)\longrightarrow\cdots\cdots
\longrightarrow(n,c_n,c^n_1,\dots,c^n_b,s_n)\longrightarrow v\ .$$

Intuitively, traveling on the path, being at vertex $(k,c_k,c^k_1,\dots,c^k_b,s_k)$ indicates
that the prefix $(c_1,\dots,c_k)$ of the sequence had been determined, it sums to $s_k$, and
its $b$ largest values are $c^k_1\geq\cdots\geq c^k_b$ and satisfy the majorization
inequalities $\sum_{j=1}^i c^k_j\leq\sum_{j=1}^i d_j$ for $i=1,\dots,b$.
Formally, by definition of the edges in $D$, we have that the tuple $c=(c_1,\dots,c_n)$
corresponding to the path satisfies $\sum_{j=1}^n c_j=s_n=\sum_{j=1}^n d_j$, and,
permuting it to be nonincreasing as $c_{\pi(1)}\geq\cdots\geq c_{\pi(n)}$,
it satisfies that $c_{\pi(j)}=c^n_j$ for $j=1,\dots,b$ and hence
$\sum_{j=1}^i c_{\pi(j)}\leq\sum_{j=1}^i d_j$ for $i=1,\dots,b$.
Since the row sums are bounded by $b$, as explained in the beginning of the proof of
Part 2 above, it follows that $c_{\pi(1)}\geq\cdots\geq c_{\pi(n)}$ is majorized by $d$.
So by Proposition \ref{Ryser}, there is a matrix $A\in\{0,1\}^{m\times n}$ with row sum
tuple $r$ and column sum tuple $c$. By the definition of the lengths, the
length of this directed path is $\sum_{j=1}^n f_j(c_j)=\sum_{j=1}^n f_j(c_j(A))$.

Conversely, for any matrix $A\in\{0,1\}^{m\times n}$ with row sum tuple $r$ and column sum
tuple $c$, its nonincreasing permutation $c_{\pi(1)}\geq\cdots\geq c_{\pi(n)}$
is majorized by $d$ by Proposition \ref{Ryser}, hence the following is a
directed path from $u$ to $v$ in $D$ of length $\sum_{j=1}^n f_j(c_j)=\sum_{j=1}^n f_j(c_j(A))$,
$$u\longrightarrow(1,c_1,c_1,0,\dots,0,c_1)\longrightarrow \cdots\cdots
\longrightarrow(n,c_n,c_{\pi(1)},\dots,c_{\pi(b)},\sum_{j=1}^n c_j)\longrightarrow v\ .$$

So the column sum problem indeed reduces to the shortest directed path problem.

\vskip.2cm
Now, the number of vertices of $D$ is bounded by $n(m+1)^{b+1}(mn+1)+2=O(m^{b+2}n^2)$ and since
$b$ is fixed it is polynomial in the input. So a shortest directed path from $u$ to $v$
in $D$ can then be obtained in polynomial time, see e.g. \cite{Sch}. Then an optimal
nonincreasing tuple $c=(c_1,\dots,c_n)$ of column sums can be read off from this optimal
directed path as explained above. Finally, an optimal matrix $A\in\{0,1\}^{m\times n}$,
with row sum tuple $r=(r_1,\dots,r_m)$ and column sum tuple $c=(c_1,\dots,c_n)$,
can be computed in polynomial time by Lemma \ref{flow}.
\epr

\section*{Acknowledgments}
S. Onn was supported by a grant from the Israel Science Foundation and by the Dresner chair
at the Technion. He thanks the referees for useful comments that improved the presentation.

\end{document}